\documentclass[amstags,fleqn,reqno]{article}

\usepackage{diagrams}

\usepackage{amsthm}
\usepackage{amsmath}
\usepackage{amssymb}

\theoremstyle{plain}
\newtheorem{theorem}{Theorem}[section]
\newtheorem{lemma}[theorem]{Lemma}
\newtheorem{corollary}[theorem]{Corollary}

\theoremstyle{definition}
\newtheorem{definition}[theorem]{Definition}

\newtheorem{construction}[theorem]{Construction}

\newtheorem{example}[theorem]{Example}

\usepackage{xspace}

\renewcommand\geq{\geqslant}
\renewcommand\leq{\leqslant}

\newcommand{\fpf}{fixed-point-free\xspace} 
\newcommand{\fp}{fixed-point\xspace}

\newcommand{\Tidentity}{\textnormal{(T1)}\xspace}
\newcommand{\Tcycledisjoint}{\textnormal{(T2)}\xspace}
\newcommand{\Tfpf}{\textnormal{(T3)}\xspace}
\newcommand{\Ttransitive}{\textnormal{(T4)}\xspace}

\newcommand{\Rdisjoint}{\textnormal{(R1)}\xspace}
\newcommand{\RtripA}{\textnormal{(R2)}\xspace}
\newcommand{\RtripB}{\textnormal{(R3)}\xspace}

\newcommand{\darts}{\Omega\xspace} 
\newcommand{\settriangles}{\mathcal{T}}

\newcommand{\triToT}{\psi}

\newcommand{\opa}{\diamond} 
\newcommand{\opb}{\otimes}

\newcommand{\blackvertex}{\bullet} 
\newcommand{\whitevertex}{\circ} 
\newcommand{\starvertex}{\star} 

\newcommand{\transversal}{\mathfrak{T}}

\newcommand{\mov}{\textnormal{Mov}}

\usepackage{ifpdf}

\ifpdf
        \usepackage[pdftex]{graphicx}
\else
	\usepackage{epsfig}
        \usepackage{graphicx}
\fi
\title{Partitioning $3$-homogeneous latin bitrades}

\author{Carlo H\"{a}m\"{a}l\"{a}inen \\
Department of Mathematics\\ The University of Queensland\\
St Lucia, 4072 \\
Brisbane, Australia\\
Telephone: +61 7 3346 1431\\ 
Facsimile: +61 7 3365 1677 \\
carloh@maths.uq.edu.au}

\begin{document}

\maketitle

\begin{abstract}
A latin bitrade $(T^{\diamond},\, T^{\otimes})$ 
is a pair of partial latin squares that define the difference between
two arbitrary latin squares $L^{\diamond} \supseteq T^{\diamond}$
and $L^{\otimes} \supseteq T^{\otimes}$ of the same order. A $3$-homogeneous bitrade 
$(T^{\diamond},\, T^{\otimes})$ has three
entries in each row, three entries in each column, and each symbol
appears three times in $T^{\diamond}$. Cavenagh~\cite{MR2241536} showed that any
$3$-homogeneous bitrade may be partitioned into three transversals.
In this paper we provide an independent
proof of Cavenagh's result using geometric methods. In doing so we provide a
framework for studying bitrades as tessellations in spherical, euclidean
or hyperbolic space. Additionally, we show how latin bitrades are
related to finite representations of certain triangle groups.
% \keywords{latin square \and latin bitrade \and triangle group \and tessellation}
% \subclass{05B15 \and 05B45}
\end{abstract}

\section{Introduction}

A {\em latin bitrade} $(T^{\opa},\, T^{\opb})$ is a pair of partial
latin squares which are disjoint, occupy the same set of non-empty
cells, and whose corresponding rows and columns contain the same set
of entries. One of the earliest studies of latin bitrades appeared 
in~\cite{DrKe1}, where they are referred to as {\em exchangeable partial
groupoids}. Latin bitrades are prominent in the study of 
{\em critical sets}, which are minimal defining sets of latin squares
(\cite{BaRe2},\cite{codose},\cite{Ke2}) and the intersections
between latin squares (\cite{Fu}). We write $i \opa j = k$ when
symbol $k$ appears in the cell at the intersection of row $i$ and column
$j$ of the (partial) latin square $T^{\opa}$. A $3$-homogeneous bitrade
has $3$ elements in each row, $3$ elements in each column, and each
symbol appears $3$ times. Cavenagh~\cite{MR2241536} obtained the
following theorem, using combinatorial methods,
as a corollary to a general classification result
on $3$-homogeneous bitrades. 

\begin{theorem}[Cavenagh~\cite{MR2241536}]\label{theorem3homTransversals}
Let $(T^{\opa},\, T^{\opb})$ be a $3$-homogeneous bitrade. Then
$T^{\opa}$ can be partitioned into three transversals.  
\end{theorem}

In this paper we provide an independent
and geometric proof of Cavenagh's result. In doing so we provide a
framework for studying bitrades as tessellations in spherical, euclidean
or hyperbolic space.
In particular, bitrades can be thought of as finite representations of
certain triangle groups.

We let permutations act on the right, in accordance with computer
algebra systems such as Sage~\cite{sage}.
Graphs in this paper may contain loops or multiple edges; otherwise our
notation is standard and we refer the reader to
Diestel~\cite{MR2159259}. Some basic topological terms will be used; for
these we refer the reader to Stillwell~\cite{MR1211642}. Finally, a good
reference for hypermaps and graphs on surfaces is~\cite{LZ}.

\section{Latin bitrades}

A {\em partial latin square} $P$ of order $n > 0$ is an 
$n \times n$ array where each $e \in \{ 0, 1, \dots, n-1 \}$
appears at most once in each row, and at most once in each column.
A {\em latin square} $L$ of order $n > 0$ is an
$n \times n$ array where each $e \in \{ 0, 1, \dots, n-1 \}$
appears exactly once in each row, and exactly once in each column.
It is convenient to use setwise notation to refer to entries
of a (partial) latin square, and we
write $(i,j,k) \in P$ if and only if symbol $k$ appears in the
intersection of row $i$ and column $j$ of $P$.
In this manner, $P \subseteq A_1 \times A_2 \times A_3$ for finite sets
$A_i$, each of size $n$.
It is also convenient to interpret a (partial) latin square as a multiplication
table for a binary operator $\opa$, writing
$i \opa j = k$ if and only if $(i,j,k) \in T = T^{\opa}$.

% A {\em partial latin square} $T^{\opa}$ is a set $T$ along with a binary
% operation that satisfies the following condition: if $x \opa y$ is
% defined, then $x \opa y = z$ has a unique solution if any two of $x$,
% $y$, and $z$ are fixed. It is convenient to use the setwise notation
% $(x,\, y,\, z) \in T^{\opa}$ to denote $x \opa y = z$.
% We index the rows, columns, and symbols by the finite sets
% $A_1$, $A_2$, and $A_3$, respectively.

\begin{definition}\label{defnBitradeA123}
Let $T^{\opa}$, $T^{\opb} \subseteq A_1 \times A_2 \times A_3$ be
two partial latin squares. Then $(T^{\opa},\, T^{\opb})$ is a {\em
bitrade} if the following three conditions are satisfied:
\begin{itemize}

\item[\Rdisjoint] $T^{\opa} \cap T^{\opb} = \emptyset$;

\item[\RtripA] for all $(a_1,\, a_2,\, a_3) \in T^{\opa}$ and all $r$,
$s \in \{1,\, 2,\, 3\}$,
$r \neq s$, there exists a unique $(b_1,\, b_2,\, b_3) \in T^{\opb}$
such that $a_r=b_r$ and $a_s=b_s$;

\item[\RtripB] for all $(a_1,\, a_2,\, a_3) \in T^{\opb}$ and all $r$,
$s \in \{1,\, 2,\, 3\}$,
$r \neq s$, there exists a unique $(b_1,\, b_2,\, b_3) \in T^{\opa}$
such that $a_r=b_r$ and $a_s=b_s$.

\end{itemize}
\end{definition}

Conditions~\RtripA and \RtripB imply that each row (column) of
$T^{\opa}$ contains the same subset of $A_3$ as the corresponding row
(column) of $T^{\opb}$.  A {\em $k$-homogeneous bitrade} $(T^{\opa},\, T^{\opb})$ 
has $k$ entries in each row of $T^{\opa}$, $k$ entries in each column of
$T^{\opa}$, and each symbol appears $k$ times in $T^{\opa}$. By symmetry
the same holds for $T^{\opb}$. A set $\transversal \subseteq T^{\opa}$ is a {\em
transversal\/} if $\transversal$ intersects each row of $T^{\opa}$ in precisely one
entry, each column in precisely one entry, and if the number of symbols
appearing in $\transversal$ is equal to $\left| \transversal \right|$. The latter condition
can be written as
$\left| \{ k \mid (i,\, j,\, k) \in \transversal \} \right| = \left|
\transversal \right|$.
A bitrade $(T^{\opa},\, T^{\opb})$ is {\em primary\/}
if whenever $(U^{\opa},\,U^{\opb})$ is a bitrade such that
$U^{\opa}\subseteq
T^{\opa}$ 
and
$U^{\opb}\subseteq   
T^{\opb}$, then 
$(T^{\opa},\, T^{\opb})=  
(U^{\opa},\, U^{\opb})$.  
Bijections $A_i \rightarrow A'_i$, for $i = 1$, $2$, $3$, give an 
{\em isotopic} bitrade, and permuting each $A_i$ gives an {\em autotopism}.

In~\cite{Dr9}, Dr\a'apal gave a representation of bitrades in terms of
three permutations $\tau_i$ acting on a finite set.  
For $r \in \{1,\, 2,\, 3\}$, define the map 
$\beta_r \colon T^{\opb} \rightarrow
T^{\opa}$ where $(a_1,\, a_2,\, a_3) \beta_r = (b_1,\, b_2,\, b_3)$ if and only if
$a_r \neq b_r$
and $a_i = b_i$ for $i \neq r$.
By Definition~\ref{defnBitradeA123} each $\beta_r$ is a bijection.
Then
$\tau_1,\, \tau_2,\, \tau_3\colon T^{\opa} \rightarrow
T^{\opa}$ are defined by
\begin{align}
\tau_1 &= \beta_2^{-1}\beta_3, \qquad
\tau_2 = \beta_3^{-1}\beta_1, \qquad
\tau_3 = \beta_1^{-1}\beta_2. \label{eqnTau}
\end{align}
We refer to
$[\tau_1,\, \tau_2,\, \tau_3]$
as the $\tau_i$ {\em representation}.
We write $\mov(\pi)$ for the set of points that the (finite) permutation
$\pi$ acts on.

\begin{definition}\label{defnT1234}
Let $\tau_1$, $\tau_2$, $\tau_3$ be (finite) permutations
and let $\darts = \mov(\tau_1) \cup \mov(\tau_2) \cup \mov(\tau_3)$.  
Define four properties:
\begin{enumerate}

\item[\Tidentity] $\tau_1 \tau_2 \tau_3 = 1$;

%\item[\Tcycledisjoint] an orbit $\rho$ of $\tau_i$ has at most one point in common
%with an orbit $\mu$ of $\tau_j$, where $1 \leq i < j \leq 3$;
\item[\Tcycledisjoint] if $\rho_i$ is a cycle of $\tau_i$
and $\rho_j$ is a cycle of $\tau_j$
then $\left| \mov(\rho_i) \cap \mov(\rho_j) \right| \leq 1$,
for any $1 \leq i < j \leq 3$;

\item[\Tfpf] each $\tau_i$ is \fpf;

\item[\Ttransitive] the group $\langle \tau_1,\, \tau_2,\, \tau_3 \rangle$ is
transitive on $\darts$.
\end{enumerate}

\end{definition}

By letting $A_i$ be the set of cycles of $\tau_i$, Dr\a'apal obtained
the following theorem, which relates
Definition~\ref{defnBitradeA123}
and~\ref{defnT1234}.  

\begin{theorem}[{Dr\a'apal~\cite{Dr9}}]\label{theoremDrapalTauStructure}
A bitrade $(T^{\opa},\, T^{\opb})$ is equivalent (up to isotopism) to three
permutations $\tau_1$, $\tau_2$, $\tau_3$ acting on a set $\darts$
satisfying \Tidentity, \Tcycledisjoint, and \Tfpf. 
If \Ttransitive is also satisfied then the bitrade is primary.
\end{theorem}

To construct the $\tau_i$ representation for a bitrade we simply evaluate
Equation~\eqref{eqnTau}. In the reverse direction we have the following
construction:

\begin{construction}[$\tau_i$ to bitrade]\label{constructionTauToBitrade}
Let $\tau_1$, $\tau_2$, $\tau_3$ be permutations satisfying
Condition~\Tidentity, \Tcycledisjoint, and \Tfpf.
Let $\darts = \mov(\tau_1) \cup \mov(\tau_2) \cup \mov(\tau_3)$.  
Define $A_i = \{ \rho \mid \textnormal{$\rho$ is a cycle of $\tau_i$} \}$
for $i = 1$, $2$, $3$. Now define two arrays
$T^{\opa}$, $T^{\opb}$: 
\begin{align*}
T^{\opa}= \{( {\rho}_1,\, {\rho}_2,\, {\rho}_3 )\mid 
&\text{ $\rho_i \in A_i$ and 
$\left| \mov(\rho_1) \cap \mov(\rho_2) \cap \mov(\rho_3) \right| \geq
1$} \} \\
T^{\opb} = \{( {\rho}_1,\, {\rho}_2,\, {\rho}_3 ) \mid
&\text{ $\rho_i \in A_i$ and $x$, $x'$, $x''$ are distinct points of
$\darts$ such} \\
&\text{ that $x\rho_1=x'$, $x'\rho_2=x''$, $x''\rho_3=x$} \}. 
\end{align*}
By Theorem~\ref{theoremDrapalTauStructure} $(T^{\opa},\, T^{\opb})$ 
is a bitrade.
\end{construction}

\begin{example}\label{exampleIntercalateRep}
The smallest bitrade $(T^{\opa},\, T^{\opb})$ is the {\em intercalate},
which has four entries. The bitrade is shown below:
\begin{align*}
T^{\opa} = \begin{array}{c|cccc}
\opa & 0 & 1 \\
\hline 
0 & 0 & 1 \\
1 & 1 & 0
\end{array}
& \quad & 
T^{\opb} = \begin{array}{c|cccc}
\opb & 0 & 1 \\
\hline 
0 & 1 & 0 \\
1 & 0 & 1
\end{array}
%& \quad & 
%\begin{array}{c|cccc}
%\opa & ~ & ~ \\
%\hline 
%~ & 000 & 011 \\
%~ & 101 & 110
%\end{array}
\end{align*}
The $\tau_i$ representation is
$\tau_1 = (000,011)(101,110)$, $\tau_2 = (000,101)(011,110)$,
$\tau_3 = (000,110)(011,101)$, 
where we have written $ijk$ for $(i,j,k) \in T^{\opa}$ to make
the presentation of the $\tau_i$ permutations clearer.
%%%%%%%%%%%%%%%%%%%%%%%%%%%%%%%%%%%%%%%%%%%%%
By Construction~\ref{constructionTauToBitrade}
with $\darts = \{ 000, 011, 101, 110 \}$
we can convert the $\tau_i$ representation to a bitrade 
$(U^{\opa},\, U^{\opb})$:
\begin{align*}
U^{\opa} &= \begin{array}{c|cccc}
\opa & (000,101) & (011,110) \\
\hline 
(000,011) & (000,110) & (011,101) \\
(101,110) & (011,101) & (000,110)
\end{array} \\
U^{\opb} &= \begin{array}{c|cccc}
\opb & (000,101) & (011,110) \\
\hline 
(000,011) & (011,101) & (000,110) \\
(101,110) & (000,110) & (011,101)
\end{array}
\end{align*}
In this way we see that row $0$ of $T^{\opa}$ corresponds to 
row $(000,011)$ of $U^{\opa}$, which is the cycle
$( 000,011 )$ of $\tau_1$, and so on for the columns and symbols.
\end{example}
\begin{example}\label{ex3hom}
The following $3$-homogeneous bitrade is pertinent to the proof of the
main result of this paper:
%%%%%%%%%%%%%%%%%%%%%%%%%%%%%%%%%%%%%%%%%%%%%%%%%%%%%%%%
\begin{align}
T^{\opa} = \begin{array}{c|cccc}
\opa & 1 & 2 & 3 & 4 \\
\hline 
1 & 1 & 3 & \phantom{-} &2 \\
2 & 3 & 2 & 4 &\phantom{-} \\
3 & \phantom{-} & 4 & 3 &1 \\
4 & 2 & \phantom{-} & 1 &4
\end{array}
& \quad & 
T^{\opb} = \begin{array}{c|cccc}
\opb & 1 & 2 & 3 & 4 \\
\hline 
1 & 3 & 2 &\phantom{-} &1 \\
2 & 2 & 4 &3 &\phantom{-} \\
3 & \phantom{-} & 3 &1 &4 \\
4 & 1 & \phantom{-} &4 &2
\end{array}
\label{eqn3homBitrade}
\end{align}
Writing $ijk \in T^{\opa}$ for
$(i,\, j,\, k) \in T^{\opa}$, the $\tau_i$ representation is
\begin{alignat*}{3}
\tau_1 &=(111,\, 142,\, 123) (213,\, 234,\, 222) (324,\, 341,\, 333) (412,\, 444,\, 431) \\
\tau_2 &=(111,\, 213,\, 412) (123,\, 222,\, 324) (234,\, 333,\, 431) (142,\, 341,\, 444) \notag \\
\tau_3 &=(111,\, 431,\, 341)(123,\, 333,\, 213)(142,\, 412,\,
222)(234,\, 444,\, 324) \notag  \\
\darts &= \mov(\tau_1) \cup \mov(\tau_2) \cup \mov(\tau_3).\notag 
\end{alignat*}
%
% Calculating tau3 in GAP:
%
% tau1 := (111, 142, 123) (213, 234, 222) (324, 341, 333) (412, 444, 431);;
% tau2 := (111, 213, 412) (123, 222, 324) (234, 333, 431) (142, 341, 444);;
% tau3 := Inverse(tau1*tau2);;
%
The bitrade has four rows so $\tau_1$ has four cycles; 
similarly $\tau_2$ and $\tau_3$ each have four cycles. (In general,
a bitrade can have a different number of row, column, and
symbol cycles.)
Using Construction~\ref{constructionTauToBitrade},
the cell at row $(111,\, 142,\, 123)$,
column $(111,\, 213,\, 412)$, will contain the symbol
$(111,\, 431,\, 341)$ since these cycles intersect in $111$. 
\end{example}

\section{Bitrades as graphs on surfaces}

% I suggest that in the first paragraph of Section 3 you add construction
% 1.4.16 from page 14 of your thesis.
% 
% Then probably consturction 1.4.26 also expand on your Example 3.1
% to show exactly how you get your bipartite graph and then your
% triangulation from that.

Before showing how a bitrade can be represented as a graph embedded in a
surface, we briefly review the theory of hypermaps.  A {\em combinatorial hypermap\/} 
$[\sigma,\, \alpha,\, \varphi]$
is made up of three permutations
$\sigma$, $\alpha$, $\varphi$ and
a finite set $\darts$ such that
$\sigma \alpha \varphi = 1$ and $G = \langle \sigma, \alpha
\rangle$ acts transitively on $\darts$. The following construction takes
a combinatorial hypermap to a {\em hypermap}, which is a 
bipartite graph embedded in a surface. For
a proof of correctness see Chapter~1 of~\cite{LZ} and references therein, 
and for further examples see Chapter~1 and~2 of~\cite{hamalainen2007}.
The representation of hypermaps as bipartite graphs was given by
Walsh~\cite{MR0360328}.

\begin{construction}\label{constrPairToBipartiteGraph}
Let $[\sigma,\, \alpha,\, \varphi]$ be a combinatorial hypermap on the
finite set~$\darts$.
Create vertex sets $V_1$, $V_2$ and undirected edges~$E$:
\begin{align*}
V_1 &= \{ v \mid \text{$v$ is an cycle of $\sigma$} \} \\
V_2 &= \{ v \mid \text{$v$ is an cycle of $\alpha$} \} \\
E &= \{ \{ v,v' \}_x \mid \text{$v \in V_1$, $v' \in V_2$, and 
$x \in \mov(v) \cap \mov(v')$ }\}
\end{align*}
Colour the vertices of $V_1$ black (denoted~$\blackvertex$) and those of
$V_2$ white (denoted~$\whitevertex$). When drawing the graph we usually 
label an edge $\{ v,v' \}_x$ with $x$ to save space.
Suppose that $(x_1, x_2, \dots, x_n)$ is a cycle of $\sigma$, and let $v$ be
the associated black vertex with adjacent edges
$\{ v, v_i \}_{x_i}$ for $1 \leq i \leq n$.
Then order the edges adjacent to $v$ as 
$\{ v, v_1 \}_{x_1}$,
$\{ v, v_2 \}_{x_2}$,
\dots,
$\{ v, v_n \}_{x_n}$ in the anticlockwise direction. Apply the same
process to each $v \in V_2$. This defines a {\em rotation scheme} for the
vertices of the bipartite graph, and hence an embedding in a surface.
\end{construction}

\begin{example}\label{exBipartiteGraph}
Let $\darts = \{ a,b,c,d \}$ and define
$\sigma = (a,b)(c,d)$ and
$\alpha = (a,c)(b,d)$.
Then there are two black vertices, two white vertices, and four edges:
%\begin{align*}
$V_1 = \{ (a,b),\, (c,d) \}$
$V_2 = \{ (a,c),\, (b,d) \}$ and
$E = \{ a, b, c, d \}$.
%\end{align*}
The graph embedding, with anticlockwise orientation, is shown 
in Figure~\ref{figBipartiteEmbedding}.
%%%%%%%%%%%%%%%%%%%%%%%%%%%%%%%%%%%%%%%%%%%%%%%%%%
\begin{figure}
\begin{center}
\includegraphics{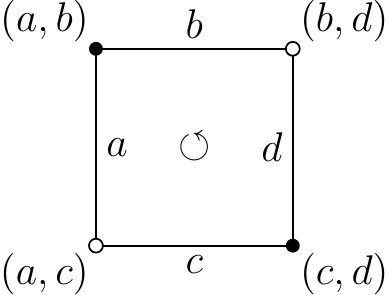}
\caption{An embedding of a bipartite graph.}\label{figBipartiteEmbedding}
\end{center}
\end{figure}
%%%%%%%%%%%%%%%%%%%%%%%%%%%%%%%%%%%%%%%%%%%%%%%%%% 
\end{example}

Given a bipartite graph embedding, we often move to the {\em canonical
triangulation}, as described in the following construction:

%%%%%%%%%%%%%%%%%%%%%%%%%%%%%%%%%%%%%%%%%%%%%%%%%%%%%%%%%%%%%%%%%%%%
\begin{construction}[{\cite[p.~50]{LZ}}]\label{constrFundamentalTriangulation}
Let $\mathcal{H}$ be a hypermap.
Place a new vertex~$\starvertex$ in each face of the hypermap. Connect this
new vertex to each vertex that lies on the border of the face
using
dotted edges to~$\blackvertex$ vertices
and
dashed edges to~$\whitevertex$ vertices. The
surface is now subdivided into triangles. Each triangle has three types of
vertices: $\blackvertex$, $\whitevertex$, and $\starvertex$; each
triangle has three types of
sides: a solid, dashed, or dotted line.
From the inside of a triangle, we view
its vertices according to the 
order $\blackvertex$,~$\whitevertex$,~$\starvertex$,~$\blackvertex$, and
if we turn in the
anticlockwise
direction then the triangle is {\em
positive}, otherwise it is {\em negative}. We shade the positive
triangles.
\end{construction}

%%%%%%%%%%%%%%%%%%%%%%%%%%%%%%%%%%%%%%%%%%%%%%%%%%
\begin{figure}
\begin{center}
    \includegraphics{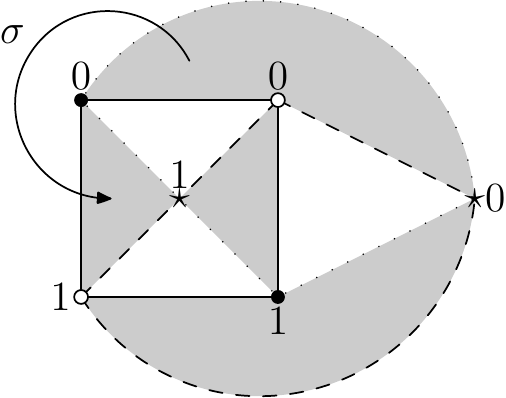}
\caption{Canonical triangulation of the bipartite graph embedding of
Figure~\ref{figBipartiteEmbedding}. Vertices are labelled from the
set $\{0,1 \}$ to aid in identification with the intercalate bitrade,
and the action of $\sigma$ is shown on a particular shaded triangle.}\label{figCanonical}
\end{center}
\end{figure}
\begin{figure}
\begin{center}
% {\includegraphics{glueHypermap-3.pdf}}
{\includegraphics{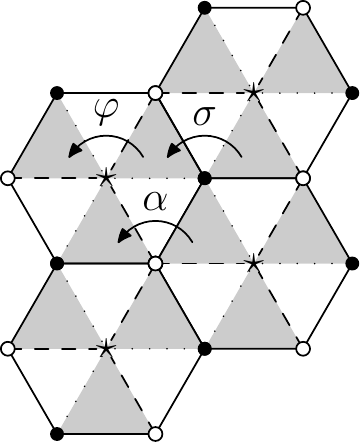}}
\caption{Canonical triangulation of a particular hypermap.  The action
of $\sigma$, $\alpha$, and $\varphi$ on certain shaded triangles is
shown.}\label{figFundTriangulation}
\end{center}
\end{figure}

%%%%%%%%%%%%%%%%%%%%%%%%%%%%%%%%%%%%%%%%%%%%%%%%%% 
Since each (shaded) triangle is adjacent to precisely one solid edge in
the canonical triangulation, we can identify the action of $\sigma$ as
the rotation of shaded triangles around their black vertex in an
anticlockwise direction, as shown in 
Figure~\ref{figCanonical} (also, see \cite[p.~51]{LZ}). In general,
the action of
$\alpha$ and $\varphi$ correspond to rotations around
white and star vertices, as indicated in 
Figure~\ref{figFundTriangulation}.
%%%%%%%%%%%%%%%%%%%%%%%%%%%%%%%%%%%%%%%%%%%%%%%%%%%%%%%%%%%%%%%%%%%%
%%%%%%%%%%%%%%%%%%%%%%%%%%%%%%%%%%%%%%%%%%%%%%%%%%%%%%%%%%%%%%%%%%%%
\begin{example}
The canonical triangulation of the bipartite graph embedding
of Example~\ref{exBipartiteGraph} is shown 
in Figure~\ref{figCanonical}.
Writing $ijk$ for the shaded triangle with vertex labels $i$, $j$, $k$
on black, white and star vertices, respectively, we see that
$000 \sigma  = 011$,
$000 \alpha  = 101$, and
$000 \varphi = 110$.
As expected, the action of $\sigma$, $\alpha$, and $\varphi$ in 
Figure~\ref{figCanonical}
is exactly the same
as $\tau_1$, $\tau_2$, and $\tau_3$ of
Example~\ref{exampleIntercalateRep}.  
\end{example}
%%%%%%%%%%%%%%%%%%%%%%%%%%%%%%%%%%%%%%%%%%%%%%%%%%%%%%%%%%%%%%%%%%%%
Applying Euler's formula leads to the {\em genus formula\/} for
hypermaps:
\begin{equation}\label{eqnHypermapGenus}
z(\sigma)+z(\alpha)+z(\varphi) - \left| \Omega \right| = 2 - 2g
\end{equation}
where $z(\pi)$ denotes the number of cycle of the permutation $\pi$.

\begin{lemma}\label{lemmaTorus}
A $3$-homogeneous bitrade $(T^{\opa},\, T^{\opb})$ defines a
tessellation of shaded and unshaded triangles in the Euclidean plane.
Each shaded triangle is edge-wise adjacent only to unshaded triangles
(and vice-versa). Shaded and unshaded triangles correspond to
the entries of $T^{\opa}$ and $T^{\opb}$, respectively. Black, white,
and star vertices correspond to row, column, and symbol labels of
$T^{\opa}$.
\end{lemma}

\begin{proof}
Let $(T^{\opa},\, T^{\opb})$ be a $3$-homogeneous bitrade
and let $[\tau_1,\, \tau_2,\, \tau_3]$ be the $\tau_i$ representation.
By Condition~\Tidentity and ~\Ttransitive we see that
$\tau_1$, $\tau_2$, and $\tau_3$ satisfy the properties to be a
combinatorial hypermap.
Let $[\sigma,\, \alpha,\, \varphi] = [\tau_1,\, \tau_2,\, \tau_3]$ and
construct the associated hypermap using
Construction~\ref{constrPairToBipartiteGraph}.
Apply Construction~\ref{constrFundamentalTriangulation} so that the
hypermap consists of shaded and unshaded triangles.  
Since $z(\tau_1) = z(\tau_2) = z(\tau_3) = \left| T^{\opa} \right|/3$
and $\left| \Omega  \right| = \left| T^{\opa} \right|$ it follows that
$g=1$ so the underlying surface is the torus.
The fundamental group of the torus is 
$\mathbb{Z} \times \mathbb{Z}$ so the covering surface
is the Euclidean plane.
By Construction~\ref{constrFundamentalTriangulation},
each shaded triangle is adjacent, edge-wise, to precisely one unshaded
triangle, and vice-versa. The permutation $\sigma$ acts on
shaded triangles while $\tau_1$ acts on elements of $T^{\opa}$
by Equation~\eqref{eqnTau}. We set $\sigma = \tau_1$ so
shaded triangles correspond to elements of $T^{\opa}$
and unshaded triangles correspond to elements of $T^{\opb}$.
Black vertices correspond to cycles of $\tau_1$ which, in turn, correspond
to row labels of $T^{\opa}$ (and similar for white and star vertices).
\end{proof}

\begin{figure}
\begin{center}
\includegraphics{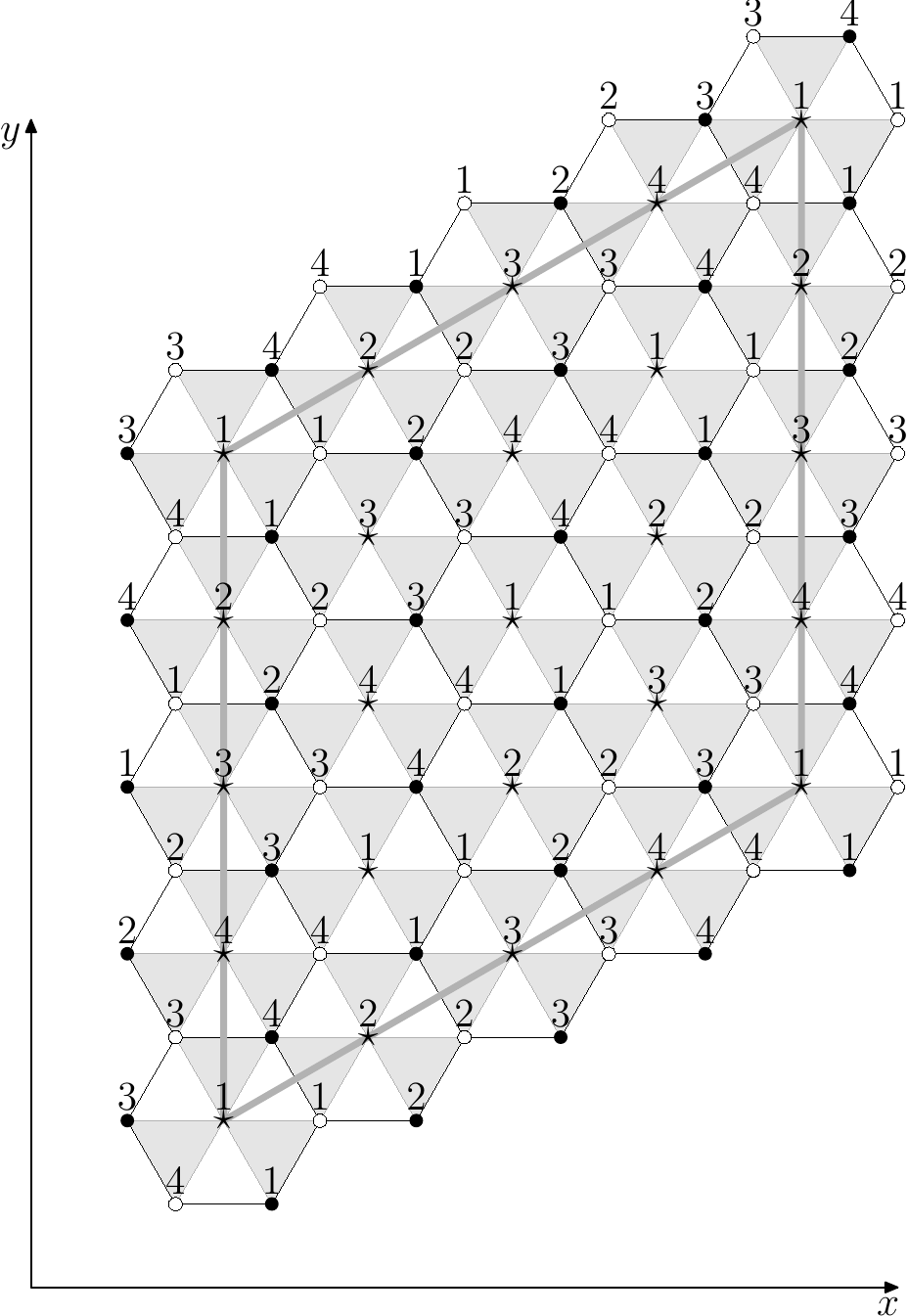}
\end{center}
\caption{Canonical triangulation of the bitrade in
Example~\ref{ex3hom}.
Identifying opposite sides of the solid grey parallelogram gives a
torus. \label{fig3homTessellation}}
\end{figure}

\begin{example}
Figure~\ref{fig3homTessellation} shows the tessellation for the
$3$-homogeneous bitrade of Example~\ref{ex3hom}.  Identifying opposite
sides of the parallelogram marked by thick grey lines gives the torus.
With regards to Theorem~\ref{theorem3homTransversals},
we can partition $T^{\opa}$ into three transversals
\begin{align*}
\transversal_1 &= \{ (1,1,1), (2,2,2), (3,3,3), (4,4,4) \}, \\
\transversal_2 &= \{ (1,4,2), (2,1,3), (3,2,4), (4,3,1) \}, \\
\transversal_3 &= \{ (1,2,3), (2,3,4), (3,4,1), (4,1,2) \}
\end{align*}
where $T^{\opa} = \transversal_1 \cup \transversal_2 \cup \transversal_3$.
These transversals may be located geometrically in
Figure~\ref{fig3homTessellation}: 
$\transversal_1$ is made up of shaded triangles located directly above a $\blackvertex$ vertex,
$\transversal_2$ is made up of shaded triangles located directly to the
lower-left of a $\blackvertex$ vertex,
and 
$\transversal_3$ is made up of shaded triangles located directly to the
lower-right of a $\blackvertex$ vertex.
\end{example}
%%%%%%%%%%%%%%%%%%%%%%%%%%%%%%%%%%%%%%%%%%%%%%%%%%%%%
%%%%%%%%%%%%%%%%%%%%%%%%%%%%%%%%%%%%%%%%%%%%%%%%%%%%%

\section{The geometric proof}
\label{section3homGeometricProof}

In this section we provide the geometric proof of
Theorem~\ref{theorem3homTransversals}. Let $(T^{\opa},\, T^{\opb})$ be a
$3$-homogeneous bitrade. Apply Lemma~\ref{lemmaTorus} to obtain the
labelled tessellation of the Euclidean plane for $(T^{\opa},\, T^{\opb})$.
Without loss of generality, let the triangles have unit length sides.
%Let $t$ be one such an
%equilateral triangle with unit length sides lying in the Euclidean
%plane. Mark its three vertices by $\blackvertex$, $\whitevertex$, and
%$\starvertex$.
%Give the Euclidean plane an anticlockwise orientation. 
%Shade the triangle, since under the given orientation,
%vertices are visited in the order
%$\blackvertex$,~$\whitevertex$,~$\starvertex$,~$\blackvertex$.
Let $t$ be an unshaded triangle in the tessellation.
Define three actions $\rho_i$ on $t$:
$\rho_1$ rotates $t$ by angle $2 \pi/3$ anticlockwise around its
$\blackvertex$ vertex;
$\rho_2$ rotates $t$ by angle $2 \pi/3$ anticlockwise around its
$\whitevertex$ vertex;
$\rho_3$ rotates $t$ by angle $2 \pi/3$ anticlockwise around its
$\starvertex$ vertex.
%%%%%%%%%%%%%%%%%%%%%%%%%%%
%% The action of the $\rho_i$ leaves three unshaded triangles around each vertex. Extend the
%% action of $\rho_i$ as follows, where $t'$ is an unshaded triangle:
%% $\rho_1$ rotates $t'$ by angle $2 \pi/3$ clockwise around its
%% $\blackvertex$ vertex;
%% $\rho_2$ rotates $t'$ by angle $2 \pi/3$ clockwise around its
%% $\whitevertex$ vertex;
%% $\rho_3$ rotates $t'$ by angle $2 \pi/3$ clockwise around its
%% $\starvertex$ vertex.
%% The $\rho_i$ rotate unshaded triangles in a clockwise direction just as
%% $\sigma$, $\alpha$, and $\varphi$ of a hypermap rotate unshaded triangles
%% in the clockwise direction.
%%%%%%%%%%%%%%%%%%%%%%%%%%%
The plane is tessellated by hexagons similar to those in 
Figure~\ref{fig3homTessellation}.

Since $(T^{\opa},\, T^{\opb})$ is $3$-homogeneous, it follows that there
are three shaded triangles at each vertex, so
$\rho_i^3 = 1$ for $1 \leq i \leq 3$ and
$\rho_1 \rho_2 \rho_3 = 1$. These $\rho_i$ induce a 
{\em triangle group\/} $\Gamma$ which acts on the
set of equilateral triangles $\settriangles$ of the tessellation:
\begin{equation}\label{eqnTessellationGamma}
\Gamma = \langle \rho_1, \rho_2, \rho_3 \mid \rho_1^3 = \rho_2^3 = \rho_3^3 
= \rho_1 \rho_2 \rho_3 = 1 \rangle.
\end{equation}
If $[\tau_1, \tau_2, \tau_3]$ is the $\tau_i$ representation for the
bitrade in question, then we define the {\em cartographic group} $G$ by:
\begin{equation}\label{eqnTessellationG}
G = \langle \tau_1, \tau_2, \tau_3 \mid \tau_1^3 = \tau_2^3 = \tau_3^3 
= \tau_1 \tau_2 \tau_3 = \ldots = 1 \rangle.
\end{equation}
Note that $\Gamma$ is an infinite group, acting on the tessellation of
the Euclidean plane, while 
$G$ is a finite permutation group, acting on the corresponding
triangles on the identified surface (the torus). 
The group $G$ has all of the defining relations for $\Gamma$ so it is
natural to define a group homomorphism $\theta$ that sends $\rho_i$ to $\tau_i$ and the empty word 
$1_{\Gamma}$ to the identity $1_G$. We then extend
$\theta$ to an arbitrary word by
$\rho_{i_1} \rho_{i_2} \cdots \rho_{i_n} \mapsto 
\tau_{i_1} \tau_{i_2} \cdots \tau_{i_n}$
where $i_{\ell} \in \{ 1,2,3 \}$ for $1 \leq \ell \leq n$.

To relate the group actions 
$\settriangles \times \Gamma \rightarrow \settriangles$
and 
$T^{\opa} \times G \rightarrow \darts$
we form a map $\triToT \colon \settriangles \rightarrow \darts$.
Fix a shaded triangle $t_0 \in \settriangles$ and 
an entry $x_0 \in \darts$ and set $t_0 \triToT = x_0$. Then use $\theta$
to extend $\triToT$ to any $t \in \settriangles$ by
defining
$t\triToT = x_0 (\delta \theta)$ where $t_0 \delta = t$ for some $\delta \in \Gamma$.

\begin{lemma}\label{3homCommutes}
The map $\triToT$ with base points $t_0$ and $x_0$ is well defined and
commutes with the actions of $\Gamma$ and $G$.
\end{lemma}

\begin{proof}
Let $\triToT$ be defined as above for some fixed 
$t_0$, $x_0$. 
First we check that $\triToT$ is a well-defined map, namely that the
choice of $\delta$ for $t_0 \delta = t$ does not matter.
Suppose that $t_0 \delta_1 = t = t_0 \delta_2$. Then
$t_0 \delta_1 \delta_2^{-1} = t_0$ so 
$(\delta_1 \delta_2^{-1}) \theta = g \in G$ where
$x_0g = x_0$ (we can't assume that $g$ is the identity in $G$, only that
it fixes $x_0$). Then
\begin{align*}
&x_0 g = x_0 (\delta_1 \delta_2^{-1}) \theta = x_0 (\delta_1 \theta) (\delta_2^{-1} \theta) = x_0  \\
\Rightarrow \quad & x_0 (\delta_1 \theta) = x_0 (\delta_2^{-1}
\theta)^{-1} = x_0 (\delta_2 \theta)  \\
\Rightarrow \quad & x_0 (\delta_1 \theta) = x_0 (\delta_2 \theta)
\end{align*} 
so $t \triToT$ takes the same value whether $\delta_1$ or $\delta_2$ was
chosen. Hence $\triToT$ is well defined.

Next we check that $\triToT$ commutes with both group actions
$\nu$ and $\eta$  of $\Gamma$ and $G$, respectively. In other words, the
following diagram must commute: 
%%%%%%%%%%%%%%%%%%%%%%%%%%%%%%%%%%%%%%%
	\begin{diagram}
	\settriangles \times \Gamma & \rTo^{\nu} & \settriangles \\
	\dTo^{\triToT \times \theta} &                 & \dTo_{\triToT} \\
	\darts \times G & \rTo^{\eta} & \darts 
	\end{diagram}
%%%%%%%%%%%%%%%%%%%%%%%%%%%%%%%%%%%%%%%
Choose 
$t \in \settriangles$, $\xi \in \Gamma$. Then
\begin{align*}
&(t, \xi) \nu \triToT = (t \xi) \triToT = x_0 (\delta_1 \theta) \\
&(t, \xi) (\triToT \times \theta) \eta = 
(x_0 (\delta_2 \theta), \xi \theta) \eta
= x_0 (\delta_2 \theta) (\xi \theta)
= x_0 ((\delta_2 \xi) \theta)
\end{align*}
where $t_0 \delta_1 = (t \xi)$ and $t_0 \delta_2 = t$.
Then $t_0 (\delta_2 \xi) = t \xi = t_0 \delta_1$ so
$\nu \triToT = (\triToT \times \theta) \eta$ and the diagram commutes.  
\end{proof}

\begin{figure}
\begin{minipage}{0.4\linewidth}
\centering
	\includegraphics{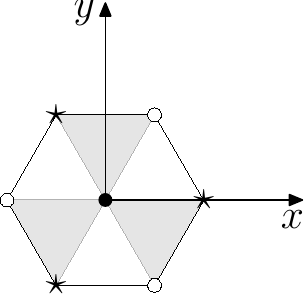}
	\caption{One possible placement of the $x$, $y$ axes. \label{fig3homXYaxes}  }
\end{minipage}
\hspace{0.5cm} % To get a little bit of space between the figures
\begin{minipage}{0.5\linewidth}
\centering
	\includegraphics{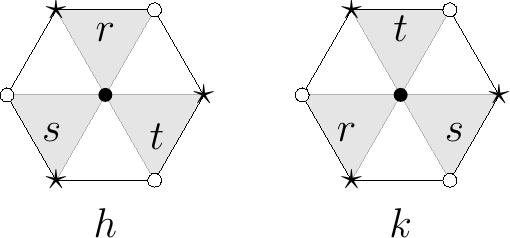}
	\caption{Inconsistently labelled hexagons due to 
	$r \in \transversal_1 \cap \transversal_2$.\label{fig3homConsistency}}
\end{minipage}
\end{figure}

The tessellation lies on the Euclidean plane and we are free to place
the $x$, $y$ axes as we wish. We will choose one of three placements:
that shown in Figure~\ref{fig3homXYaxes}, or the rotation of those axes
by angle $2 \pi/3$ or $-2 \pi/3$. 
%%%%%%%%%%%%%%%%%%%%%%%%%%%%%%%%%%%%%%%%%%%%%%%%
%%%%%%%%%%%%%%%%%%%%%%%%%%%%%%%%%%%%%%%%%%%%%%%%
In other words, the origin is always
at a $\blackvertex$ vertex, the $x$ axis always
`points' through a $\starvertex$ vertex, and the $y$ axis (in the positive
direction) bisects a shaded triangle.
\begin{definition}\label{defThreeTransversals}
Geometrically define three subsets $\transversal_i \subset T^{\opa}$ 
as follows:
\begin{equation*}
\begin{split}
\transversal_1 = \{ {t}\triToT \mid \; &\text{$t$ is the shaded triangle immediately above $v$ with $v$ as its $\blackvertex$ vertex} \} \\
\transversal_2 = \{{t}\triToT \mid  \;&\text{$t$ is the shaded triangle immediately to} \\
		&\text{the lower-left of $v$  with $v$ as its $\blackvertex$ vertex} \} \\
\transversal_3 = \{{t}\triToT \mid  \;&\text{$t$ is the shaded triangle immediately to} \\
&\text{the lower-right of $v$  with $v$ as its $\blackvertex$ vertex} \}
\end{split}
\end{equation*} 
where $v$ ranges over all~$\blackvertex$ vertices.
\end{definition}

In what follows it will be convenient to label a triangle $t$ in the
tessellation by $t \triToT$. We now move on to showing that the
$\transversal_i$ sets are mutually disjoint.

\begin{lemma}\label{3homPartition12}
Let $\transversal_1$ and $\transversal_2$ be given as in
Definition~\ref{defThreeTransversals}. Then
$\transversal_1 \cap \transversal_2 = \emptyset$.  
\end{lemma}

\begin{proof}
Suppose, for contradiction, that there exists a triple $r \in T^{\opa}$ such that
$r \in \transversal_1 \cap \transversal_2$. 
Since $\tau_1$ has no \fp (Condition~\Tfpf in the definition of
a bitrade) there must be a
$3$-cycle $(r,\, s,\, t)$ in $\tau_1$ for some $s$, $t \in T^{\opa}$. If 
$u \triToT = r$ for some shaded triangle $u \in \settriangles$
then it must be that
$(u \rho_1) \triToT = r \tau_1 = s$
and
$(u \rho_1^2) \triToT = r \tau_1^2 = t$.
Recalling that $\rho_1$ is rotation about a $\blackvertex$~vertex in the
anticlockwise direction, 
the labelled tessellation must have hexagons like those shown in 
Figure~\ref{fig3homConsistency}. 
The hexagon $h$ has $r \in \transversal_1$ while the hexagon $k$
has $r \in \transversal_2$. The order of $r$, $s$, $t$ is forced by $\triToT$.

Let each side of a triangle in the tessellation have unit length.
Without loss of
generality we can place the $x$, $y$ axes on the tessellation (as in
Figure~\ref{fig3homXYaxes}) so that the $\blackvertex$~vertex of $h$
is at $(0,0)$ and the $\blackvertex$~vertex of $k$ is at $(x,y)$.
Then we have a Euclidean distance
$d(h,k) = \sqrt{x^2 + y^2}$ which we assume to be minimal.  
We then show that there exists
another pair of inconsistently labelled hexagons $h'$ and $k'$ such that
$d(h',k') < d(h,k)$ except for a few cases in which contradictions
arise with respect to the bitrade itself. In the limiting case we get
$d(h',k') = 0$ which implies that $\tau_1$ has a fixed point $r$.
There are four main cases to check, each with three subcases
{\bf a}, {\bf b}, and {\bf c}.
Each of the {\bf a} and {\bf b} cases cover an infinite part of the plane
so we use various constructions to find $h'$, $k'$ such that
$d(h',k') < d(h,k)$. The {\bf c} cases are finite and provide the
required local contradictions.

\paragraph{Case 1: $x$, $y \geq 0$.}

\paragraph{Case 1a: $x > 3/2$, $y \geq 0$.}

Suppose that $t \tau_1 \tau_2 = w$
and consider the action of $\rho_1 \rho_2$ on the triangles labelled $t$
as shown in Figure~\ref{figcase1.1}.
Recall that $\rho_1$ is rotation to the next shaded triangle in an
anticlockwise direction around a $\blackvertex$~vertex, and 
$\rho_2$ is rotation around a $\whitevertex$~vertex.
We find $d(h',k')$ and factor out the $d(h,k)$ term:
\begin{align*}
d(h',k')^2	&= (x-3)^2 + y^2 \\
		&= x^2 - 6x + 9 + y^2 \\
		&= d(h,k)^2 - 6x + 9.
\end{align*}
Now $d(h',k')^2 - d(h,k)^2 =  - 6x + 9 < 0$
since $x > 3/2$
so
$d(h',k') < d(h,k)$. The last step is to observe that there must be a
cycle $(w,\, w',\, w'')$ in $\tau_1$ as shown above. Then
$w'' \in \transversal_1 \cap \transversal_2$, so $h'$ and $k'$ are a
closer pair of inconsistent hexagons. This completes~Case~1a.

\begin{figure}
%\begin{minipage}{0.4\linewidth}
%\centering
	%\includegraphics{threetrans-1.pdf}
	%\caption{Marked hexagon from Figure~\ref{fig3homTessellation}.\label{fig3homPattern}}
%\end{minipage}
%\hspace{0.2cm} % To get a little bit of space between the figures
%\begin{minipage}{0.4\linewidth}
\centering
	\includegraphics{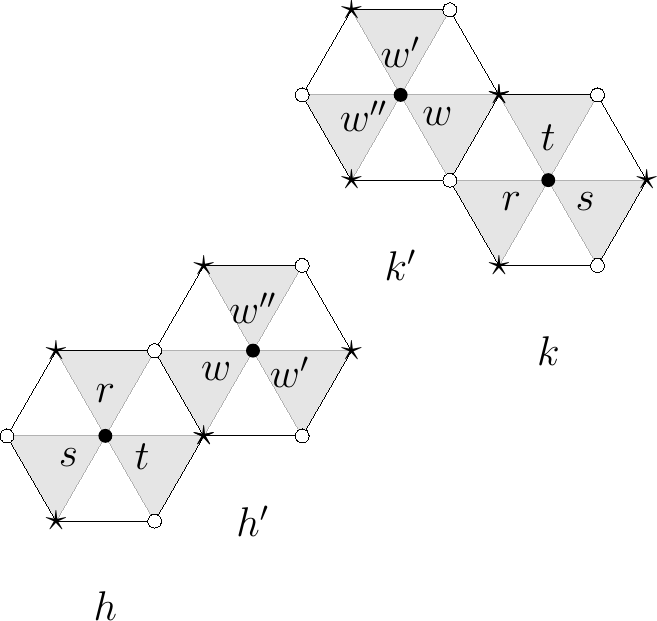}
	\caption{Case 1a.\label{figcase1.1}}
%\end{minipage}
\end{figure}
%%%%%%%%%%%%%%%%%%%%%%%%%%%%%%%%%%%%%%%%%%%
\begin{figure}
\begin{minipage}{0.4\linewidth}
\centering
	\includegraphics{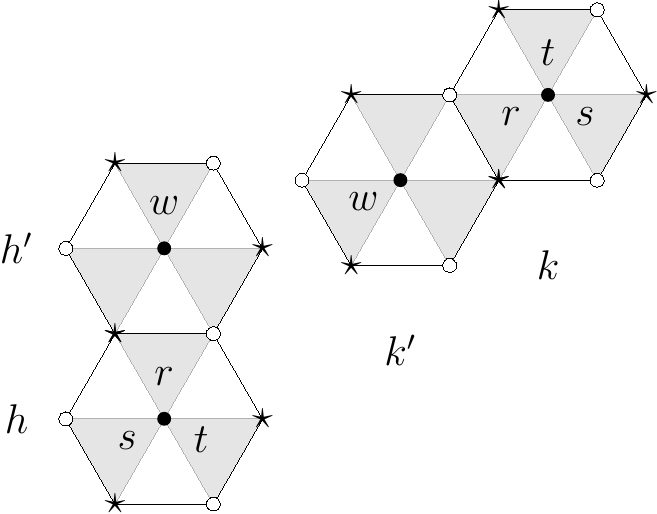}
	\caption{Case 1b.\label{figcase1.2}}
\end{minipage}
\hspace{1.3cm} % To get a little bit of space between the figures
\begin{minipage}{0.6\linewidth}
\centering
	\includegraphics{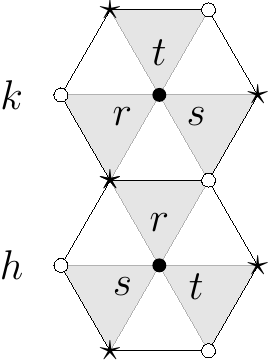}
	\caption{Case 1c, (ii).\label{figcase1.3}}
\end{minipage}
\end{figure}
%%%%%%%%%%%%%%%%%%%%%%%%%%%%%%%%%%%%%%%%%%%
\begin{figure}
\begin{minipage}{0.5\linewidth}
\centering
	\includegraphics{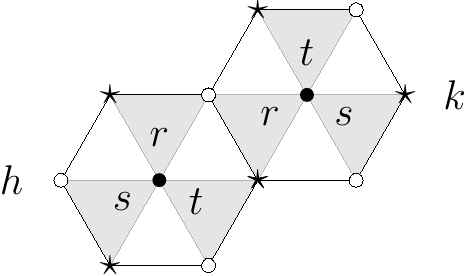}
	\caption{Case 1c, (iii).\label{figcase1.4}}
\end{minipage}
\hspace{-0.51cm} % To get a little bit of space between the figures
\begin{minipage}{0.6\linewidth}
\centering
	\includegraphics{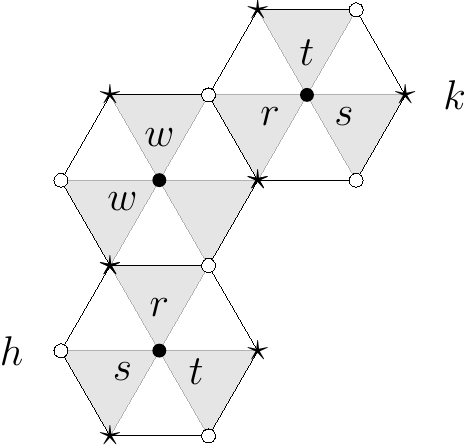}
	\caption{Case 1c, (iv).\label{figcase1.5}}
\end{minipage}
\end{figure}
%%%%%%%%%%%%%%%%%%%%%%%%%%%%%%%%%%%%%%%%%%%

% \clearpage

\begin{figure}
\begin{minipage}{0.5\linewidth}
\centering
	\includegraphics{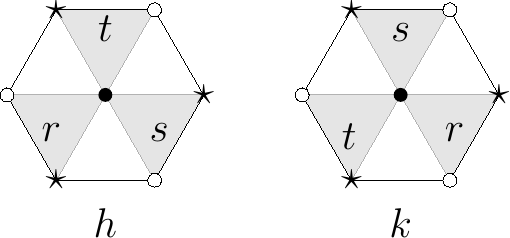}
	\caption{Inconsistently labelled hexagons due to 
		$r \in \transversal_2 \cap \transversal_3$.\label{fig3homHexagon2}}
\end{minipage}
\hspace{0.1cm} % To get a little bit of space between the figures
\begin{minipage}{0.5\linewidth}
\centering
	\includegraphics{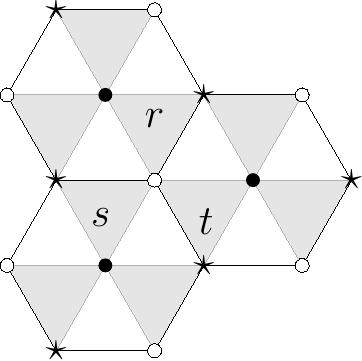}
	\caption{Consistency at a~$\whitevertex$ vertex.\label{fig3homRowCol}}
\end{minipage}
\end{figure}

% Avoids a horrendous page break in the first paragraph
% of Case 1b.
\clearpage

\paragraph{Case 1b: $0 \leq x \leq 3/2$, $y \geq 4 \frac{\sqrt{3}}{2}$.}
Suppose that $r \tau_3 \tau_1^{-1} = w$ and consider the action of 
$\rho_3 \rho_1^{-1}$ on the triangles labelled $r$ as in
Figure~\ref{figcase1.2}.
We calculate the distance between $h'$ and $k'$:
\begin{align*}
d(h',k')^2	&= (x-3/2)^2 + \left( y-3 \frac{\sqrt{3}}{2} \right)^2 \\
		&= x^2 - 3x + 9/4 + y^2 - 3\sqrt{3}y + 27/4 \\
		&= d(h,k)^2 - 3x + 9 - 3\sqrt{3}y
\end{align*}
We need $-3x + 9 - 3\sqrt{3}y < 0$ which simplifies to
$-x + 3 - \sqrt{3}y < 0$. By assumption,
$-x + 3 - \sqrt{3}y \leq -x + 3 - 4\sqrt{3}\frac{\sqrt{3}}{2} 
\leq -x + 3 - 6 = -x - 3 < 0$ and this completes Case~1b.

\paragraph{Case 1c: $0 \leq x \leq 3/2$, $0 \leq y < 4 \frac{\sqrt{3}}{2}$.}

Here we deal with local cases which give rise to contradictions. Note
that some $(x,y) \in \{ 0, \frac{3}{2} \} \times 
\{ 0, 1 \frac{\sqrt{3}}{2}, 2 \frac{\sqrt{3}}{2}, 3 \frac{\sqrt{3}}{2} \}$
do not correspond to valid hexagon positions, e.g. there is no hexagon
centred at $(0, 1 \frac{\sqrt{3}}{2})$. The valid cases are as follows:

\begin{description}

\item[(i) $x=0$, $y=0$:] In this case $h$ and $k$ are the same hexagon, implying
that $r \tau_1 = r$ so $\tau_1$ has a \fp which
contradicts~\Tfpf in the definition of a bitrade

\item[(ii) $x=0$, $y=2 \frac{\sqrt{3}}{2}$:]
In this case $\tau_3$ has a \fp $r$, contradicting~\Tfpf (see
Figure~\ref{figcase1.3}).

\item[(iii) $x=3/2$, $y=\frac{\sqrt{3}}{2}$:]
Here we find a \fp of $\tau_2$, contradicting~\Tfpf (see
Figure~\ref{figcase1.4}).

\item[(iv) $x=3/2$, $y=3\frac{\sqrt{3}}{2}$:]
If $r \tau_3 \tau_1^{-1} = w$ then the action of 
$\rho_3 \rho_1^{-1}$ on the triangles labelled $r$ shows that $\tau_1$
has a \fp as shown in Figure~\ref{figcase1.5}.

\end{description}

Cases~2, 3, and 4 are very similar. For each sub-case of type
{\bf a} and {\bf b} we state the 
$\rho_i \rho_j$ word along with the corresponding action as
$\tau_i \tau_j$.  Each sub-case of type {\bf c} gives an immediate
contradiction in the form of a fixed point for
some $\tau_1$, $\tau_2$,  or $\tau_3$.

{\bf Case 2:} $x$, $y \leq 0$.

{\bf Case 2a:} $x < -3/2$, $y \leq 0$.
Use $\rho_3 \rho_1^{-1}$ where $s \tau_3 \tau_1^{-1} = w$.

{\bf Case 2b:} $-3/2 \leq x \leq 0$, $y < -3 {\sqrt{3}}/{2}$.
Use $\rho_1 \rho_2$ where $s \tau_1 \tau_2 = w$.

{\bf Case 2c:} $-3/2 \leq x \leq 0$, $-3 {\sqrt{3}}/{2} \leq y \leq 0$.

{\bf Case 3:} $x \geq 0$, $y \leq 0$.

{\bf Case 3a:} $x > 3/2$, $y \leq 0$.
Use $\rho_2 \rho_1^{-1}$ where $r \tau_2 \tau_1^{-1} = w$.

{\bf Case 3b:} $0 \leq x \leq 3/2$, $y < -3 {\sqrt{3}}/{2}$.
Use $\rho_3 \rho_1^{-1}$ where $t \tau_3 \tau_1^{-1} = w$.

{\bf Case 3c:} $0 \leq x \leq 3/2$, $-3 {\sqrt{3}}/{2} \leq y \leq 0$.

{\bf Case 4:} $x \leq 0$, $y \geq 0$.

{\bf Case 4a:} $x < -{3}/{2}$, $y \geq 0$.
Use $\rho_2^{-1} \rho_1$, where $s \tau_2^{-1} \tau_1 = w$.

{\bf Case 4b:} $-{3}/{2} \leq x \leq 0$, $y > 2 {\sqrt{3}}/{2}$.
Use $\rho_2 \rho_1^{-1}$, where $s \tau_2 \tau_1^{-1} = w$.

{\bf Case 4c:} $-{3}/{2} \leq x \leq 0$, $0 \leq y \leq 2
{\sqrt{3}}/{2}$. \\

This completes the proof of Lemma~\ref{3homPartition12} 
\end{proof}

\begin{corollary}\label{3homPartition123}
Let $\transversal_1$, $\transversal_2$, and $\transversal_3$, be given as in
Definition~\ref{defThreeTransversals}. Then the $\transversal_i$ are
mutually disjoint.
\end{corollary}

\begin{proof}
Suppose that there exists $r \in T^{\opa}$ such that 
$r \in \transversal_2 \cap \transversal_3$. Since $\tau_1$ is \fpf it
must contain a $3$-cycle $(r,\, s,\, t)$ for some $s$, $t \in T^{\opa}$.
Then the inconsistent
hexagons are as shown in Figure~\ref{fig3homHexagon2}. 
We see that $t \in \transversal_1 \cap \transversal_2$, so by
Lemma~\ref{3homPartition12} we have a contradiction. The case where
$r \in \transversal_1 \cap \transversal_3$ is similar.
\end{proof}

\begin{corollary}\label{corPartition}
Let $\transversal_1$, $\transversal_2$, and $\transversal_3$, be given as in
Definition~\ref{defThreeTransversals}. Then the set $\{ \transversal_1,
\transversal_2, \transversal_3 \}$ is a partition of $T^{\opa}$.  
\end{corollary}

\begin{proof}
By Corollary~\ref{3homPartition123} the $\transversal_i$ sets are
mutually disjoint.
By assumption, the group $G = \langle \tau_1, \tau_2, \tau_3 \rangle$ acts transitively
on $T^{\opa}$. The group homomorphism $\theta$ given by
$\rho_i \mapsto \tau_i$ is actually a group epimorphism. With Lemma~\ref{3homCommutes}
it follows that each $r \in T^{\opa}$ will be an element of some
$\transversal_i$ set. Hence 
$\{ \transversal_1, \transversal_2, \transversal_3 \}$ is a partition of
$T^{\opa}$.
\end{proof}

\begin{lemma}\label{3homTransversal}
The three sets $\transversal_1$, $\transversal_2$, $\transversal_3$ as defined in
Definition~\ref{defThreeTransversals} are transversals.
\end{lemma}

\begin{proof}
In
Lemma~\ref{3homPartition12} and
Corollary~\ref{corPartition} we constructed the partition
by dividing up elements of $T^{\opa}$ around a $\blackvertex$~vertex, so
it is impossible for any
$\transversal_i$ to have more than one element from a row of $T^{\opa}$
(in particular this would imply a fixed point of $\tau_1$).
Conversely, suppose that a cycle $(r,\, s,\, t)$ exists in $\tau_2$ and is
labelled as shown in Figure~\ref{fig3homRowCol}.
Now $r \in \transversal_3$, $s \in \transversal_1$, $t \in \transversal_2$ according to the labelling
induced around $\blackvertex$ vertices. If another hexagon centred at a
$\whitevertex$ vertex was inconsistently labelled then we would have an
inconsistent labelling around $\blackvertex$ vertices, contradicting
Corollary~\ref{corPartition}. Similarly, labellings around $\starvertex$ vertices
are consistent.
\end{proof}

Corollary~\ref{corPartition} and Lemma~\ref{3homTransversal} give the
main result, Theorem~\ref{theorem3homTransversals}.
We note that, in general,
the covering surface
will be spherical, Euclidean, or hyperbolic. Most (large) bitrades will be
hyperbolic, and we expect that future work will derive combinatorial
properties of hyperbolic bitrades from their geometrical representation.

%\bibliographystyle{plain}
%\bibliography{hamalainen}

\end{document}